\documentclass[12pt]{article}
\usepackage{amsfonts,amssymb,amsmath,amsbsy}
\usepackage{latexsym}
\usepackage{graphicx}
\usepackage{tikz}

\newtheorem{theo}{Theorem}
\newtheorem{lemma}[theo]{Lemma}
\newtheorem{prop}[theo]{Proposition}

\newtheorem{conj}[theo]{Conjecture}
\newtheorem{coro}[theo]{Corollary}

\newcommand{\proof}{\noindent{\it Proof: }}
\newcommand{\proofbox}{\hfill \mbox{ $\Box$}\\}

\newcommand{\R}{{\mathbb R}}
\newcommand{\C}{{\mathbb C}}

\title{Iterated functions and the Cantor set in one dimension}
\author{Benjamin Hoffman}

\begin{document}
\maketitle

\section{Introduction}

\begin{figure} 
   \centering
   \includegraphics[width=6in]{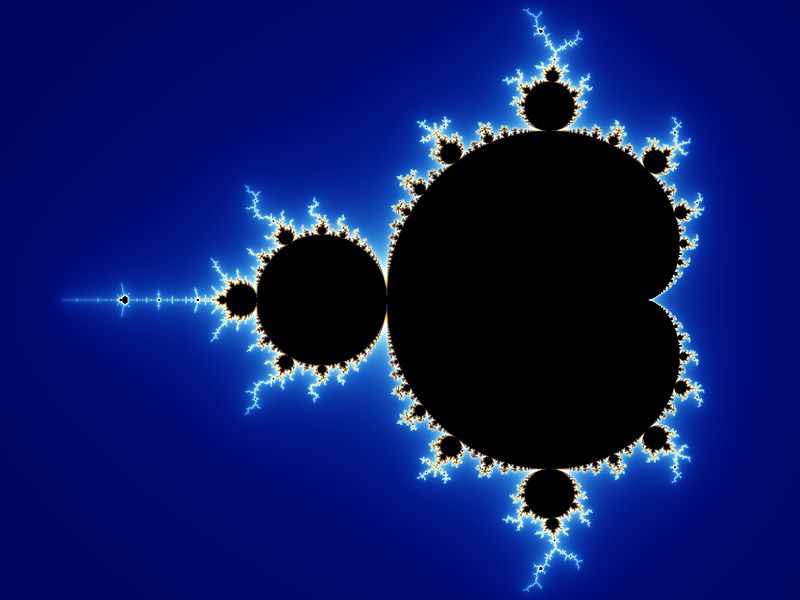} 
   \caption{The Mandelbrot set. Points within the set are black, points outside the set are colored based on how quickly $F_c$ diverges. Image created by Wolfgang Beyer and used under the Creative Commons Attribuition-Share Alike 3.0 Unported license.}
   \label{Figure0}
\end{figure}

In this paper we consider the long-term behavior of points in $\R$ under iterations of continuous functions. For some motivation as to why iterated functions are interesting to study, we can look to the \emph{Mandelbrot Set} (Figure \ref{Figure0}). Consider the function $F_c: \C \to \C$ given by $F_c(z)=z^2+c$. If we let $z=0$, let $F_c^1=F_c$, and let $F_c^{n+1}=F_c\circ F_c^n$, then the behavior of $F_c^n(0)$ as $n\to \infty$ will depend on the value of $c$. The Mandelbrot set is the set of values for $c\in\C$ such that $F_c^n(0)$ remains bounded as $n\to \infty$. In Figure \ref{Figure0}, these are the points of the complex plane that are colored black. The Mandelbrot set is the canonical example of a \emph{fractal}. As you zoom in on one area of the boundary new, finer detail will always emerge. The points outside the set are colored based on how quickly $F_c^n$ diverges as $n$ grows. The beauty and the simplicity of the Mandelbrot set helped popularize the study of fractals and iterated functions in the 1980s.  For a more in depth discussion of iterated functions see \cite{devaney}. Sections \ref{secpre} and \ref{cantor} of this paper draw from the discussion there.

In this paper, we will relate iterated functions to a different kind of set. It is common in an undergraduate analysis course to construct the \emph{Cantor middle-thirds set}, by removing a countably infinite set of successively smaller open intervals from the set $[0,1]$. The Cantor middle-thirds set is an example of a certain kind of topological space called a \emph{Cantor set}. In this paper, we show that, given any Cantor set $\Lambda^*$ embedded in $\R$, there exists a continuous function $F^*:\R\to\R$ such that the points that are bounded under iterations of $F^*$ are just those points in $\Lambda^*$. In the course of this, we find a striking similarity between the way in which we construct the Cantor middle-thirds set, and the way in which we find the points bounded under iterations of certain continuous functions.

After defining the relevant terms in Section \ref{cantor}, we further motivate our main result. We prove it in Section \ref{general}. More formally, our main result is:

\begin{theo} Given a Cantor set $\Lambda^*\subset\R$, there exists a continuous function $F^*:\R\to\R$ such that both of the following hold:

(i) if $x\in\Lambda^*$ then $(F^*)^n(x)\in\Lambda^*$ for all $n$,

(ii) if $x\not\in\Lambda^*$ then $(F^*)^n(x)$ diverges to infinity.
\end{theo}

\section{Preliminaries}
\label{secpre}

Let $F:\R\to\R$ be a continuous function, let $F^1=F$ and $F^{n+1}=F\circ F^n$. Note that, depending on what $F$ is, some points $x\in\R$ are such that $|F^n(x)|\to\infty$ as $n\to\infty$, while some points are such that $F^n(x)\in[a,b]$ for all $n$, where $a,b\in\R$. In the first case, we say that $F^n(x)$ \emph{diverges to infinity}. We choose a specific function and study this type of behavior in that function.

\begin{figure}
\begin{center}
\begin{tikzpicture}[scale=.7]
\draw [help lines, <->] (-4.5,0) -- (4.5,0);
\draw [help lines, <->] (0,-4.5) -- (0,4.5);
\draw [domain=-4.5:4.5] plot (\x, {\x});
\draw [domain=-2:2] plot (\x, {\x*\x+.5});
\node [right] at (4.5,4.5) {$f$};
\node [left] at (-2,4) {$F_{\frac{1}{2}}$};
\draw  (0,.5) -- (.5,.5) -- (.5,.75) -- (.75,.75) -- (.75,1.0625) -- (1.0625,1.0625) -- (1.0625,1.629) -- (1.629,1.629) -- (1.629,3.153) -- (3.153,3.153) -- (3.153,4.5);
\end{tikzpicture}
\end{center}
\caption{Graphical analysis showing that $F_{\frac{1}{2}}(0)\to\infty$ as $n\to\infty$.} \label{fig1}
\end{figure}
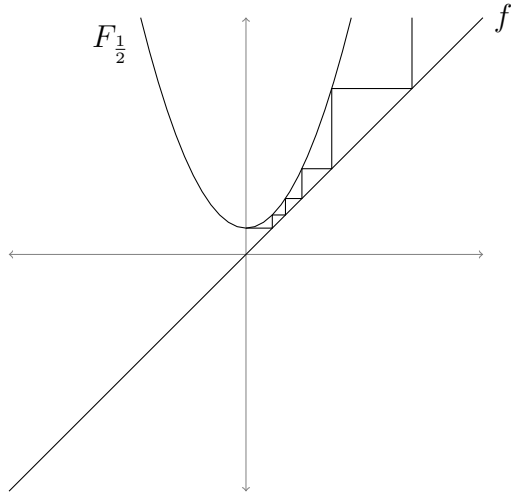

Let $F_c:\R\to\R$ be given by $F_c(x)=x^2+c$. We can study the long-term behavior of $F_c$ by plotting $F_c$ alongside the function $f(x)=x$. Given $x$, we start at $F_c(x)$, draw a line horizontally to the graph of $f$, then vertically to the graph of $F_c$, then horizontally to the graph of $f$, and so on. This heuristic method of \emph{graphical analysis} gives a quick check for the long-term behavior of points under iteration of $F_c$. For instance, we see from graphical analysis that when $c>\frac{1}{4}$, $F_c^n(x)$ diverges to infinity for all $x$ (Figure \ref{fig1}). When $c=\frac{1}{4}$, $x=\frac{1}{2}$ is a fixed point of $F_c$, since $F_c(\frac{1}{2})=\frac{1}{2}$. $F^n_c(x)$ diverges to infinity for all $x\notin[-\frac{1}{2},\frac{1}{2}]$.

\begin{figure}
\begin{center}
\begin{tikzpicture}[scale=.7]
\draw [help lines, <->] (-4.5,0) -- (4.5,0);
\draw [help lines, <->] (0,-4.5) -- (0,4.5);
\draw [domain=-4.5:4.5] plot (\x, {\x});
\draw [domain=-2.345:2.345] plot (\x, {\x*\x-1});
\node [right] at (4.5,4.5) {$f$};
\node [left] at (-2.345,4) {$F_{-1}$};
\draw (1.618,-.07) node[below]{$p$} -- (1.618,.07);
\draw (-1.618,-.07) node[below]{$-p$} -- (-1.618,.07);
\draw (-1.75,2.0625) -- (2.0625,2.0625) -- (2.0625,3.254) -- (3.254,3.254) -- (3.254,4.5);
\end{tikzpicture}
\end{center}
\caption{Graphical analysis showing that for $x\notin[-p,p]$, $F_{-1}(x)$ diverges to infinity.} \label{fig2}
\end{figure}
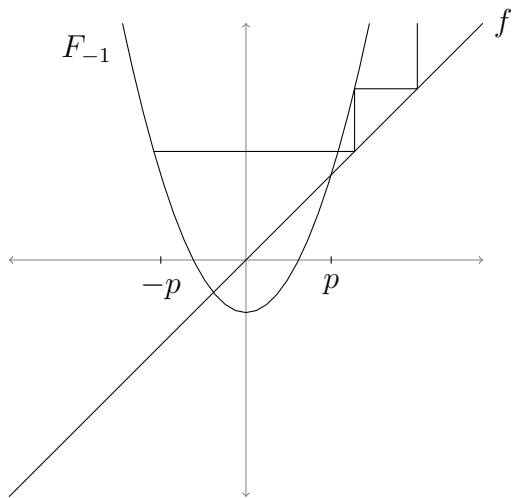

When $c<\frac{1}{4}$, note that there are two distinct solutions to $F_c(x)=x$. Let $p$ denote the larger of the two. Graphical analysis tells us that, letting $I=[-p,p]$, if $x\not\in I$, then $F^n_c(x)$ diverges to infinity (Figure \ref{fig2}).

\begin{figure}
\begin{center}
\begin{tikzpicture}[scale=.7]
\draw [help lines, <->] (-4.5,0) -- (4.5,0);
\draw [help lines, <->] (0,-4.5) -- (0,4.5);
\draw [domain=-4.5:4.5] plot (\x, {\x});
\draw [domain=-2.74:2.74] plot (\x, {\x*\x-3});
\node [right] at (4.5,4.5) {$f$};
\node [left] at (-2.74,4) {$F_{-3}$};
\draw [dashed] (2.303,0) node[below right]{$p$} -- (2.303,2.303) -- (-2.303,2.303) -- (-2.303,0) node[below left]{$-p$} -- (-2.303,-2.303) -- (2.303,-2.303) -- (2.303,0);
\draw (-.75,-2.44) -- (-2.44,-2.44) -- (-2.44,2.9536) -- (2.9536,2.9536) -- (2.9536,4.5);
\draw (-.835,0) -- node[above left]{$A_0$} (.835,0);
\draw [fill=white] (-.835,0) circle (0.1);
\draw [fill=white] (.835,0) circle (0.1);
\end{tikzpicture}
\end{center}
\caption{Graphical analysis showing that there exists an open interval $A_0\subset[-p,p]$ such that for $x\in A_0$, $F_{-3}(x)$ diverges to infinity.} \label{fig3}
\end{figure}
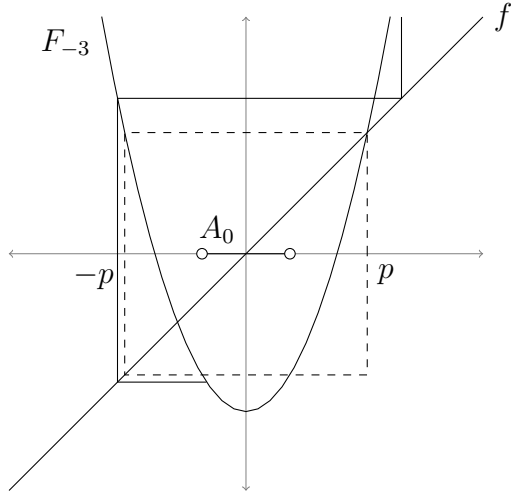

Note that, when $c<-2$, the set $A_0=\{x\in I: F_c(x)<-p\}$ is nonempty (Figure \ref{fig3}), and that if $x\in A_0$ then $F_c^n(x)$ diverges to infinity. Graphically, we find $A_0$ by drawing a square centered at the origin with sides parallel to the $x$ and $y$-axis, and which passes through the points $(p,0)$ and $(-p,0)$. $A_0$ is just those points $x\in I$ such that $F_c(x)$ lies \emph{outside} this square. Note that $A_0$ is an open interval, and that $F_c([-p,0]\backslash A_0)=F_c([0,p]\backslash A_0)=I$. Now, for $n>0$ we define $A_n=\{x\in I: F_c^n(x)\in A_0\}$.

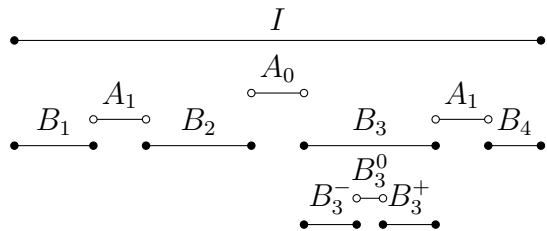
\begin{figure}
\begin{center}
\begin{tikzpicture}[scale=.7]
\draw (0,0) -- node[above]{$I$} (10,0);
\draw (4.5,-1) -- node[above]{$A_0$} (5.5,-1);
\draw (0,-2) -- node[above]{$B_1$} (1.5,-2);
\draw (1.5,-1.5) -- node[above]{$A_1$} (2.5,-1.5);
\draw (2.5,-2) -- node[above]{$B_2$} (4.5,-2);
\draw (5.5,-2) -- node[above]{$B_3$} (8,-2);
\draw (8,-1.5) -- node[above]{$A_1$} (9,-1.5);
\draw (9,-2) -- node[above]{$B_4$} (10,-2);
\draw (5.5,-3.5) -- node[above]{$B_3^-$} (6.5,-3.5);
\draw (6.5,-3) -- node[above]{$B_3^0$} (7,-3);
\draw (7,-3.5) -- node[above]{$B_3^+$} (8,-3.5);
\draw [fill=black] (0,0) circle (0.07);
\draw [fill=black] (10,0) circle (0.07);
\draw [fill=black] (0,-2) circle (0.07);
\draw [fill=black] (1.5,-2) circle (0.07);
\draw [fill=black] (2.5,-2) circle (0.07);
\draw [fill=black] (4.5,-2) circle (0.07);
\draw [fill=black] (5.5,-2) circle (0.07);
\draw [fill=black] (8,-2) circle (0.07);
\draw [fill=black] (9,-2) circle (0.07);
\draw [fill=black] (10,-2) circle (0.07);
\draw [fill=black] (6.5,-3.5) circle (0.07);
\draw [fill=black] (5.5,-3.5) circle (0.07);
\draw [fill=black] (7,-3.5) circle (0.07);
\draw [fill=black] (8,-3.5) circle (0.07);
\draw [fill=white] (1.5,-1.5) circle (0.07);
\draw [fill=white] (2.5,-1.5) circle (0.07);
\draw [fill=white] (4.5,-1) circle (0.07);
\draw [fill=white] (5.5,-1) circle (0.07);
\draw [fill=white] (8,-1.5) circle (0.07);
\draw [fill=white] (9,-1.5) circle (0.07);
\draw [fill=white] (6.5,-3) circle (0.07);
\draw [fill=white] (7,-3) circle (0.07);
\end{tikzpicture}
\end{center}
\caption{An illustration of the argument made in Proposition \ref{2nint}, for $n=1$.} \label{fig4}
\end{figure}

\begin{prop} \label{2nint} $A_n$ consists of $2^n$ open intervals (Figure \ref{fig4}). \end{prop}
\proof Assume for induction that, for some $n$, $A_n$ consists of $2^n$ open intervals and $I\backslash \bigcup_{m=0}^n A_m$ consists of $2^{n+1}$ closed intervals $B_j$, $1\le j\le 2^{n+1}$ such that for each $j$, $F^n_c(B_j)=F_c([-p,0]\backslash A_0)$ or $F^n_c(B_j)=F_c([0,p]\backslash A_0)$. Then we have that $B_j$ consists of two closed intervals $B_j^-,B_j^+$ and an open interval $B_j^0$ such that $F^{n+1}_c(B_j^-)=[-p,0]\backslash A_0$, $F^{n+1}_c(B_j^+)=[0,p]\backslash A_0$, and $F^{n+1}_c(B_j^0)=A_0$. Then $A_{n+1}$ consists of the open intervals $B_j^0$, $1\le j\le2^{n+1}$ and $I\backslash \bigcup_{m=0}^{n+1} A_m$ consists of $2^{n+2}$ closed intervals $B_j^-,B_j^+$, $1\le j\le 2^{n+1}$, such that each of these intervals is mapped by $F_c^{n+1}$ to $[-p,0]\backslash A_0$ or $[0,p]\backslash A_0$.\proofbox

Now, if $x\in I$ is such that there is no $n$ where $F^n_c(x)\in A_0$, then there is no point at which $x$ is mapped outside of $I$. The set $\Lambda_c=\bigcap_{n=0}^{\infty}(I\backslash A_n)=\{ x\in I:F^n_c(x)\in I {\rm{~for~all~}}n\ge 0\}$ is then the set of $x\in\R$ such that $F_c^n(x)$ does \emph{not} diverge to infinity. The construction of $\Lambda_c$ is reminiscent of the construction of  the Cantor middle-thirds set. We find that this is indeed the case; $\Lambda_c$ is a Cantor set.

\section{A Cantor set as the set of points bounded under iterations of a function} \label{cantor}

In this section, we show that for a certain function $F$, the set of points which are bounded under iterations of $F$ form a Cantor set. The process we use to find these points bears a strong resemblance to the classic construction of the Cantor middle-thirds set. A topological space $\Lambda$ is a \emph{Cantor set} iff it is totally disconnected, perfect, and compact. A topological space $X$ is \emph{totally disconnected} iff the connected components of $X$ are just the one-point sets of $X$, and \emph{perfect} iff $X$ is closed and every point in $X$ is a limit point.

\begin{prop} All Cantor sets are homeomorphic. If $\Lambda$ is a Cantor set, then $|\Lambda |$=$\mathfrak{c}$, the cardinality of the continuum. \end{prop} 

\proof See \cite{willard} for a proof of this proposition. \proofbox

Note that there exists a $c_*$ such that for all $c<c_*$ there is some $\lambda$ such that, $|F'_c(x)|>\lambda>1$ for any $x\in I\backslash A_0$. We prove that for such a $c$, $\Lambda_c$ is a Cantor set. This result holds for all $c<-2$, though this will not be proved here. For the remainder of this section, we fix $c<c_*<-2$ and drop the subscripts from $F$ and $\Lambda$.

\begin{theo} If $c<c_*$, then $\Lambda$ is a Cantor set. \end{theo}
\proof Since for all $n$, $A_n$ consists of $2^n$ open intervals, $\Lambda$ is the intersection of closed sets. Hence, $\Lambda$ is closed. Since $\Lambda\subset I=[-p,p]$, $\Lambda$ is bounded and hence compact.

We show $\Lambda$ is totally disconnected. To show that a subspace of $\R$ is totally disconnected, it suffices to show that it contains no open intervals. Assume that $[x,y]\subset\Lambda$. Then $F^n([x,y])\subset I\backslash A_0$ for all $n$. Also, by assumption, for all $z\in I\backslash A_0$, $|F'(z)|>\lambda>1$. So $|(F^n)'(z)|=|(F^{n-1})'(z)\cdot F'(F^{n-1}(z))|>\lambda|(F^{n-1})'(z)|$. Repeating this argument gives that $|(F^n)'(z)|>\lambda^n$. By the mean value theorem, there exists some $z\in[x,y]$ such that $|F^n(y)-F^n(x)|=|(y-x)(F^n)'(z)|>(y-x)\lambda^n$. $F^n$ then expands $[x,y]$ by a factor of at least $\lambda^n$, so if $x\ne y$ then there exists an $n$ such that $F^n([x,y])\not\subset I$. Thus, $x=y$ and the only closed intervals of $\Lambda$ are single point sets. Thus, $\Lambda$ contains no open intervals and is totally disconnected.

To show that $\Lambda$ is perfect, we show that if $x\in\Lambda$, then $x$ is a limit point. Consider an open interval $T$ containing $x$. By the argument in the previous paragraph, there exists an $n$ such that $-p\in F^n(T)$ or $p\in F^n(T)$. Since $-p\in\Lambda$ and $p\in\Lambda$, we have that $(F^n)^{-1}(-p)\subset \Lambda$ and $(F^n)^{-1}(p)\subset\Lambda$. Also, $((F^n)^{-1}(-p)\cup(F^n)^{-1}(p))\cap T\ne \emptyset$. Thus, any open set containing $x$ contains points in $\Lambda$ other than $x$ and hence $x$ is a limit point of $\Lambda$. Thus, $\Lambda$ is perfect. \proofbox

We now reformulate the definition of $\Lambda$ slightly, in a way that will be useful in the next section of this paper. Let $C_0=I$ and $C_{n+1}=C_n\backslash A_n$. $C_n$ then consists of $2^n$ closed intervals, and $\Lambda=\bigcap_{i=0}^n C_n$. We then obtain the following corollary:

\begin{coro} \label{CoroF} Given $\epsilon>0$, there exists $N$ such that for all $n\ge N$ the segments of $C_n$ are all less than $\epsilon$ in length.
\end{coro}

\proof Assume this is false. Then there must exist a system of intervals $\{T_n\}_{n=0}^{\infty}$, where $T_n$ is one of the $2^n$ intervals in $C_n$, and where $T_m\subset T_n$ whenever $m\ge n$. But then $\bigcap_{i=0}^{\infty}T_i\subset\Lambda$ is of length greater than $\epsilon$, contradicting the total disconnectedness of $\Lambda$. \proofbox

\section{Proof of the main theorem} \label{general}

In this section we prove the following result, generalizing some of the ideas from the last section to arbitrary Cantor sets in $\R$.

\begin{theo} \label{MainTheo} Given a Cantor set $\Lambda^*\subset\R$, there exists a continuous function $F^*:\R\to\R$ such that both of the following hold:

(i) if $x\in\Lambda^*$ then $(F^*)^n(x)\in\Lambda^*$ for all $n$,

(ii) if $x\not\in\Lambda^*$ then $(F^*)^n(x)$ diverges to infinity.
\end{theo}

We require several results for this to obtain.

\begin{prop} There exists a Cantor set $\Lambda\in\R$ and a continuous function $F$ such that (i) and (ii) in Theorem \ref{MainTheo} hold. \end{prop}

\proof $\Lambda$ and  $F$ were constructed in the previous section. \proofbox

\begin{prop} Let $\Lambda^*\subset [a,b]\in\R$ be a Cantor set such that $a,b\in\Lambda^*$. Then there is a nested sequence of closed sets $C_n^*$ such that $C_n^*$ consists of $2^n$ segments each of length less than or equal to $(\frac{2}{3})^n(b-a)$. Furthermore, $\partial(C_n^*)\subset\Lambda^*$. \end{prop} 

\proof We first show that, given $[c,d]\subset[a,b]$, $c\ne d$, there exists an open interval $(e,f)\subset[c,d]$ such that $(e,f)\cap\Lambda^*=\emptyset$. Assume not. That is, assume that there exists some $[c,d]\subset[a,b]$, $c\ne d$ such that for all open intervals $(e,f)\subset[c,d]$, $(e,f)\cap\Lambda^*\ne\emptyset$. Then there exists a point $g\in(c,d)$ such that $g\in\Lambda^*$. Because $\Lambda^*$ is totally disconnected, there is a point $h\in(c,d)$ such that $h\not\in\Lambda^*$. Then, by assumption, for all $\epsilon>0$, $(h-\epsilon,h+\epsilon)$ contains a point from $\Lambda^*$. But then $h$ is a limit point of $\Lambda^*$, and so $h\in\Lambda^*$ because $\Lambda^*$ is closed. Hence, the claim holds and given $[c,d]\subset[a,b]$, $c\ne d$, there exists an open interval $(e,f)\subset[c,d]$ such that $(e,f)\cap\Lambda^*=\emptyset$. 

Now, given $[c,d]\subset[a,b]$, choose an interval $[g,h]\subset[c,d]$ such that $g-c\ge \frac{1}{3}(d-c)$ and $d-h\ge \frac{1}{3}(d-c)$. From the previous paragraph, we may choose an open interval $(e,f)\subset[g,h]$ such that $(e,f)\cap\Lambda^*=\emptyset$, and $[c,e],[f,d]$ are such that $e-c<\frac{2}{3}(d-c)$ and $d-f<\frac{2}{3}(d-c)$. Also, given an interval $(e,f)$ that satisfies these conditions, since $\Lambda^*$ is compact, $f'=\inf\{x\in\Lambda^*:f<x\}$ and $e'=\sup\{x\in\Lambda^*:x<e\}$ exist. And $(e',f')$ also is such that $(e',f')\cap\Lambda^*=\emptyset$, $e'-c<\frac{2}{3}(d-c)$, and $d-f'<\frac{2}{3}(d-c)$. Thus, we may assume that $e,f\in\Lambda^*$.

We now define the sets $C_n^*$. Let $C_0^*=[a,b]$, which contains $2^0=1$ interval of length $(\frac{2}{3})^0(b-a)$. Given $C_n^*$ consisting of $2^n$ closed intervals of length each less than or equal to $(\frac{2}{3})^n(b-a)$ and greater than $0$, and such that $\partial(C_n^*)\subset\Lambda^*$, by the previous paragraph we may split each of these intervals into two intervals of length less than or equal to $(\frac{2}{3})^{n+1}(b-a)$. Furthermore, when we pick an open interval $(e,f)$ to remove we may do so such that $e,f\in\Lambda^*$. This guarantees that $\partial(C_{n+1}^*)\subset\Lambda^*$. Also, since $\Lambda^*$ is perfect, no intervals of $C_n^*$ are of length $0$; otherwise there would be some $x\in\Lambda^*$ that is not a limit point of $\Lambda^*$.

Let $C^*=\bigcap_{i=0}^{\infty}C_i^*$. Showing that $C^*=\Lambda^*$ completes the proof. We first fix some notation. $C_n^*$ is composed of $2^n$ closed intervals, denoted $[a_{n,j}^*,b_{n,j}^*]$ where $1\le j\le2^n$. Given $C_{n-1}^*$, to get $C_{n}^*$ we remove $2^{n-1}$ open intervals, denoted $(c_{n,j}^*,d_{n,j}^*)$, where $1\le j\le2^{n-1}$. Since there is a corresponding construction for $\Lambda$, we similarly denote the $2^n$ closed intervals of $C_n$ by $[a_{n,j},b_{n,j}]$, and the $2^{n-1}$ open intervals removed from $C_{n-1}$ by $(c_{n,j},d_{n,j})$.

Let $x\in\Lambda^*$. Then if $x\not\in C^*$, it was removed in some step of the construction of $C^*$. But this cannot be so, since we only removed points not in $\Lambda^*$. Let $x\in C^*$. Then for all $n$ there exists an interval $[a_{n,j}^*,b_{n,j}^*]\subset C_n^*$ such that $x\in[a_{n,j}^*,b_{n,j}^*]$. Since $a_{n,j}^*\in\Lambda^*$ and $|x-a_{n,j}^*|\le (\frac{2}{3})^n(a-b)$ for all $n$, $x$ is a limit point of $a_{n,j}$ and hence $x\in\Lambda^*$. Thus, $C^*=\Lambda^*$. \proofbox

\begin{coro} \label{CoroFStar} Given $\epsilon>0$, there exists $N$ such that for all $n\ge N$ the segments of $C_n^*$ are all less than $\epsilon$ in length.
\end{coro}

\proof Follows immediately from the construction of $C_n^*$. \proofbox

\begin{figure}
\begin{center}
\begin{tikzpicture}[scale=1]
\draw [help lines, ->] (0,0) -- (10,0) node[right]{$\R$};
\draw [help lines, ->] (0,0) -- (0,10) node[right]{$\R$};
\draw (2,-.5) node[below]{$a_{0,1}$} -- (8,-.5) node[right]{$C_0$} node[below]{$b_{0,1}$};
\fill (2,-.5) circle (0.05);
\fill (8,-.5) circle (0.05);
\draw (2,-1.5) node[below]{$a_{1,1}$} -- (4,-1.5) node[below]{$b_{1,1}$};
\draw (5,-1.5) node[below]{$a_{1,2}$} -- (8,-1.5) node[right]{$C_1$} node[below]{$b_{1,2}$};
\fill (2,-1.5) circle (0.05);
\fill (4,-1.5) circle (0.05);
\fill (5,-1.5) circle (0.05);
\fill (8,-1.5) circle (0.05);
\draw (2,-2.5) node[below left]{$a_{2,1}$} -- (2.5,-2.5) node[below]{$b_{2,1}$};
\draw (3.5,-2.5) node[below]{$a_{2,2}$} -- (4,-2.5) node[below right]{$b_{2,2}$};
\draw (5,-2.5) node[below]{$a_{2,3}$} -- (6,-2.5) node[below]{$b_{2,3}$};
\draw (7,-2.5) node[below]{$a_{2,4}$} -- (8,-2.5) node[right]{$C_2$} node[below]{$b_{2,4}$};
\fill (2,-2.5) circle (0.05);
\fill (2.5,-2.5) circle (0.05);
\fill (3.5,-2.5) circle (0.05);
\fill (4,-2.5) circle (0.05);
\fill (5,-2.5) circle (0.05);
\fill (6,-2.5) circle (0.05);
\fill (7,-2.5) circle (0.05);
\fill (8,-2.5) circle (0.05);
\draw(-.5,1) node[left]{$a_{0,1}^*$} -- (-.5,9) node[left]{$b_{0,1}^*$} node[above]{$C_0^*$};
\fill(-.5,1) circle (0.05);
\fill(-.5,9) circle (0.05);
\draw(-1.5,1) node[left]{$a_{1,1}^*$} -- (-1.5,4.5) node[left]{$b_{1,1}^*$};
\draw(-1.5,5.5) node[left]{$a_{1,2}^*$} -- (-1.5,9) node[left]{$b_{1,2}^*$} node[above]{$C_1^*$};
\fill (-1.5,1) circle (0.05);
\fill (-1.5,4.5) circle (0.05);
\fill (-1.5,5.5) circle (0.05);
\fill (-1.5,9) circle (0.05);
\draw(-2.5,1) node[left]{$a_{2,1}^*$} -- (-2.5,3) node[left]{$b_{2,1}^*$};
\draw(-2.5,3.5) node[left]{$a_{2,2}^*$} -- (-2.5,4.5) node[left]{$b_{2,2}^*$};
\draw(-2.5,5.5) node[left]{$a_{2,3}^*$} -- (-2.5,7) node[left]{$b_{2,3}^*$};
\draw(-2.5,7.5) node[left]{$a_{2,4}^*$} -- (-2.5,9) node[left]{$b_{2,4}^*$} node[above]{$C_2^*$};
\fill (-2.5,1) circle (0.05);
\fill (-2.5,3) circle (0.05);
\fill (-2.5,3.5) circle (0.05);
\fill (-2.5,4.5) circle (0.05);
\fill (-2.5,5.5) circle (0.05);
\fill (-2.5,7) circle (0.05);
\fill (-2.5,7.5) circle (0.05);
\fill (-2.5,9) circle (0.05);
\draw (1,0) -- node[below right]{$\phi_0$} (2,1);
\draw (8,9) -- node[below right]{$\phi_0$} (9,10);
\draw (4,4.5) -- node[below right]{$\phi_1$} (5,5.5);
\draw (2.5,3) -- node[below]{$\phi_2$} (3.5,3.5);
\draw (6,7) -- node[below]{$\phi_2$} (7,7.5);
\draw [dashed] (2.5,-2.5) -- (2.5,3) -- (-2.5,3);
\draw [dashed] (3.5,-2.5) -- (3.5,3.5) -- (-2.5,3.5);
\end{tikzpicture}
\end{center}
\caption{Definition of $\phi_0$, $\phi_1$, and $\phi_2$ based on the first three stages of the constructions of $\Lambda$ and $\Lambda^*$.} \label{bigfig}
\end{figure}
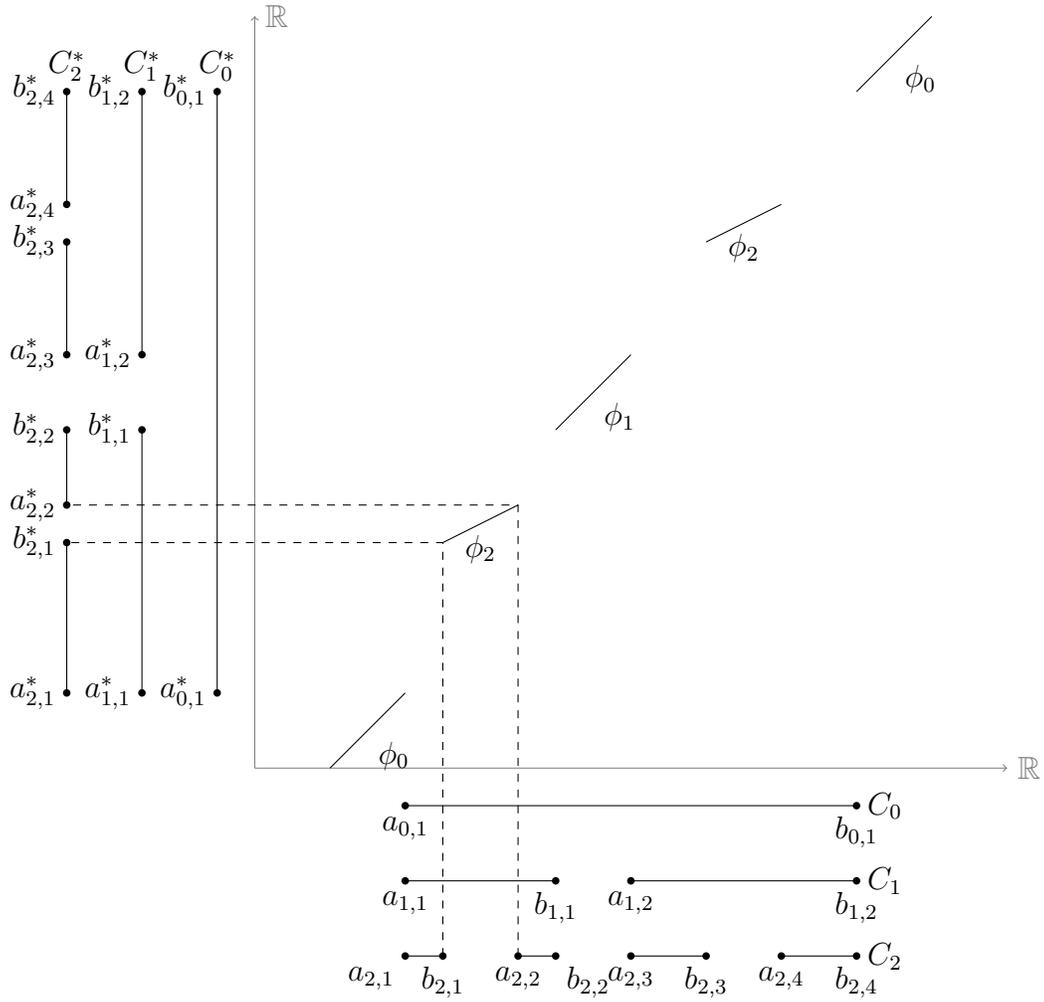

We now construct a homeomorphism $\phi:\R\to\R$ such that $\phi(\Lambda)=\Lambda^*$ (Figure \ref{bigfig}). Define $\phi_0:\R\backslash C_0\to\R\backslash C_0^*$ by
\begin{equation}
\phi_0(x) = \left\{
        \begin{array}{ll}
            x+(a_{0,1}^*-a_{0,1}) & \quad {\rm{if}}~x\le a_{0,1} \\
            x+(b_{0,1}^*-b_{0,1}) & \quad {\rm{if}}~x\ge b_{0,1}
        \end{array}
    \right.
\end{equation}
For $n\ge 1$ let $D_n={\rm{cl}}(\{(c_{n,j},d_{n,j}):1\le j\le2^n\})={\rm{cl}}(C_{n-1}\backslash C_n)$, and $D_n^*={\rm{cl}}(\{(c_{n,j}^*,d_{n,j}^*):1\le j\le2^n\})={\rm{cl}}(C_{n-1}^*\backslash C_n^*)$. Define $\phi_n:D_n\to D_n^*$ by
\begin{equation}
\phi_n(x)=\Big(\frac{d_{n,j}^*-c_{n,j}^*}{d_{n,j}-c_{n,j}}\Big)(x-c_{n,j})+c_{n,j}^* \quad {\rm{for~}}x\in[c_{i,j},d_{i,j}]{\rm{~where~}}1\le j\le2^{n-1}.
\end{equation}
Recalling that functions are sets of ordered pairs, let $\hat{\phi}=\bigcup_{n=0}^{\infty}\phi_n$. Note that $\hat{\phi}:\R\backslash\hat{\Lambda}\to\R\backslash\hat{\Lambda}^*$ is a function, where $\hat{\Lambda}=\{x\in\Lambda:\forall n,j(x\ne a_{n,j}\land x\ne b_{n,j})\}$ is the set of all points in $\Lambda$ that are not the endpoints of an interval in some $C_n$, and $\hat{\Lambda}^*=\{x\in\Lambda^*:\forall n,j(x\ne a_{n,j}^*\land x\ne b_{n,j}^*)\}$ is the analogous subset of $\Lambda^*$.

Define $\phi:\R\to\R$ by:
\begin{equation}
\phi(x) = \left\{
        \begin{array}{ll}
            \hat{\phi}(x) & {\rm{if~}} x\in\R\backslash\hat{\Lambda} \\
            \lim_{k_n\to x} \hat{\phi}(k_n) & {\rm{if~}} x\in\hat{\Lambda}
        \end{array}
    \right.
\end{equation}
where $k_n=a_{n,j}$ or $b_{n,j}$ is a fixed sequence of the boundary points of $C_n$ converging to $x$, which exists by Corollary \ref{CoroF}. The choice of $k_n$ doesn't matter; we see in Proposition \ref{ctseverywhere} that $\phi$ is continuous everywhere. The proof of this statement holds independent of our choice of $k_n$, and hence $\lim_{k_n\to x} \phi(k_n)=\lim_{k_n\to x} \hat{\phi}(k_n)$ will be the same for any sequence of points $k_n\to x$. 

We now show that $\phi$ is a homeomorphism. It suffices to show that $\phi$ is a continuous bijection.

\begin{lemma} \label{SegLemma} If $x\in[a_{n,j},b_{n,j}]\subset C_n$, then $\phi(x)\in[a_{n,j}^*,b_{n,j}^*]\subset C_n^*$. \end{lemma}

\proof If there is some $n$ such that $x\in D_n$, the result follows from the construction of $\hat{\phi}$. Otherwise, $x\in\hat{\Lambda}$. Consider the sequence $k_m\to x$ that defines $\phi(x)$. Then $k_m$ is eventually contained in the the interval $[a_{n,j},b_{n,j}]$, and since $k_m\in D_m$, $\hat{\phi}(k_m)$ is eventually contained in $[a_{n,j}^*,b_{n,j}^*]$. But $\hat{\phi}(k_m)\to\phi(x)$ by definition, and $[a_{n,j}^*,b_{n,j}^*]$ is closed, so $\phi(x)\in[a_{n,j}^*,b_{n,j}^*]$. \proofbox

\begin{prop} $\phi$ is continuous everywhere. \label{ctseverywhere} \end{prop}

\proof If $x\not\in\Lambda$, then $\phi$ is continuous on an open interval around $x$, by construction of $\hat{\phi}$. Let $x\in\Lambda$. Take a sequence $p_i\to x$. For each $p_i$, if $p_i\ne x$ let $M_i=[a_{k(i),j},b_{k(i),j}]$ be the interval of $C_{k(i)}$ containing both $x$ and $p_i$, where $k(i)$ gives the index of the last $C_n$ in which both $x$ and $p_i$ are contained in the same interval; $k(i)$ is well-defined by Corollary \ref{CoroF}. If $p_i=x$, let $M_i=\{x\}$. By Corollary \ref{CoroF}, $k(i)\to\infty$ as $p_i\to x$, whenever $k(i)$ is defined. By Lemma \ref{SegLemma}, as well as the construction of $\hat{\phi}$ when $p_i\ne x$, $\phi(p_i)\in\phi(M_i)$. By Corollary \ref{CoroFStar}, $\phi(p_i)\to\phi(x)$ and thus $\phi$ is continuous at $x$. Hence, $\phi$ is continuous everywhere. \proofbox

\begin{prop} $\phi$ is a bijection. \end{prop}

\proof Note that $\phi$ is strictly increasing on $\R\backslash\Lambda$, i.e. if $x,y\notin\Lambda$ and $x<y$, then $\phi(x)<\phi(y)$. Now, given arbitrary $x,y\in\R$, since $\Lambda$ is totally disconnected there exists sequences $x_n\to x$ and $y_n\to y$ where $x_n,y_n\notin\Lambda$ and $x_n<x$, $y_n>y$. What is more, there exists $z_1,z_2\notin\Lambda$ such that $x<z_1<z_2<y$. We then have that $\phi(x_n)<\phi(z_1)<\phi(z_2)<\phi(y_n)$ for all $n$. Hence, by continuity $\phi(x)\le\phi(z_1)<\phi(z_2)\le\phi(y_n)$. So $\phi$ is strictly increasing everywhere. Since $\phi$ is continuous, $\phi(x)\to-\infty$ as $x\to-\infty$, and $\phi(x)\to\infty$ as $x\to\infty$, it follows that $\phi$ is a bijection. \proofbox

\begin{coro} $\phi(\Lambda)=\Lambda^*$ \end{coro}

\proof This follows from the previous proposition, as well as the fact that $\phi(\R\backslash\Lambda)=\R\backslash\Lambda^*$. \proofbox

\begin{coro} $\phi:\R\to\R$ is a homeomorphism. \end{coro}

\proof Since $\phi:\R\to\R$ is a continuous bijection, by Brouwer's invariance of domain theorem (see \cite{brouwer} and page 172 of \cite{hatcher}) it is a homeomorphism. \proofbox

\begin{figure}
\begin{center}
\begin{tikzpicture}
\node (D) at (0,2.5) {$\R\backslash\Lambda\cup\Lambda$};
\node (C) at (4.5,2.5) {$\R\backslash\Lambda^*\cup\Lambda^*$};
\node (B) at (4.5,0) {$\R\backslash\Lambda^*\cup\Lambda^*$};
\node (A) at (0,0) {$\R\backslash\Lambda\cup\Lambda$};
\draw[->] (D) -- node[below]{$\phi$}  (C);
\draw[->] (D) -- node[right]{$F$} (A);
\draw[->] (A) -- node[below]{$\phi$} (B);
\draw[->] (C) -- node[right]{$F^*$} (B);
\end{tikzpicture}
\end{center}
\caption{Definition of $F^*$} \label{diagram}
\end{figure}

Now, define $F^*$ so that the diagram in Figure \ref{diagram} commutes. It is easy to check that $F^*$ is a continuous function that satisfies criteria (i) and (ii) in Theorem \ref{MainTheo}. Indeed, since $F^*=\phi\circ F\circ\phi^{-1}:\R\to\R$, $\phi$ is a homeomorphism, and $F$ is continuous, $F^*$ is continuous. If $x\in\Lambda^*$, then $\phi^{-1}(x)\in\Lambda$ by the construction of $\phi$, and so $F^n(\phi^{-1}(x))\in\Lambda$ for all $n$. Thus, $\phi(F^n(\phi^{-1}(x)))=F^*(x)\in\Lambda^*$ for all $n$, again by the construction of $\phi$. Similarly, if $x\notin\Lambda^*$, then $\phi^{-1}(x)\notin\Lambda$. For any $s\in\R$, there exists some $N$ such that for all $n>N$, $F^n(\phi^{-1}(x))>\phi^{-1}(s)$. Since $\phi$ is strictly increasing, it is order preserving and so $\phi(F^n(\phi^{-1}(x)))=(F^*)^n(x)>\phi(\phi^{-1}(s))=s$. Thus, $(F^*)^n(x)$ diverges to infinity.

We conclude by stating an additional conjecture.

\begin{conj} \label{Conj} Given a Cantor set $\Lambda^*\subset\R^m$, there exists a continuous function $F^*$ such that both of the following hold:

(i) if $x\in\Lambda^*$ then $(F^*)^n(x)\in\Lambda^*$ for all $n$,

(ii) if $x\not\in\Lambda^*$ then $||(F^*)^n(x)||$ diverges to infinity.
\end{conj}


\begin{thebibliography}{100}

\bibitem{brouwer}
{\sc Brouwer, L.}
\newblock Beweis der invarianz des $n$-dimensionalen gebiets.
\newblock {\em Mathematische Annalen}, 71 (1912), 305--315.

\bibitem{devaney}
{\sc Devaney, R.}
\newblock {\em An Introduction to Chaotic Dynamical Systems}, 2~ed.
\newblock Westview Press, 2003.

\bibitem{hatcher}
{\sc Hatcher, A.}
\newblock {\em Algebraic topology},
\newblock Cambridge University Press,
  2002.

\bibitem{willard}
{\sc Willard, S.}
\newblock {\em General Topology}.
\newblock Addison-Wesley, 1970, ch.~30.4.

\end{thebibliography}
\end{document}